\documentclass[aos,MSNbibl,nameyear,dvips]{arximspdf}

\doi{10.1214/14-AOS1205} 
\volume{42}
\issue{3}
\pubyear{2014}
\firstpage{918}
\lastpage{939}

\makeatletter

\newcommand{\rrvert}{\vert}
\newcommand{\llvert}{\vert}
\def\textonehalf{1/2}
\newtheorem{theorem}{Theorem}
\newtheorem{lemma}{Lemma}
\newtheorem{corollary}{Corollary}
\newproclaim{definition}{Definition}
\newproclaim{remark}{Remark}
\newproclaim{example}{Example}
\makeatother

\begin{document}
\begin{frontmatter}

\title{Generalized resolution for orthogonal arrays}
\runtitle{Generalized resolution for orthogonal arrays}

\begin{aug}
\author[a]{\fnms{Ulrike} \snm{Gr\"{o}mping}\corref{}\thanksref{t1}\ead[label=e1]{groemping@bht-berlin.de}\ead[label=u1,url]{http://prof.beuth-hochschule.de/groemping}}
\and
\author[b]{\fnms{Hongquan} \snm{Xu}\thanksref{t2}\ead[label=e2]{hqxu@stat.ucla.edu}\ead[label=u2,url]{http://www.stat.ucla.edu/\textasciitilde hqxu/}}
\runauthor{U. Gr\"{o}mping and H. Xu}
\thankstext{t1}{Supported in part by Deutsche
Forschungsgemeinschaft Grant GR 3843/1-1.}
\thankstext{t2}{Supported by NSF Grant DMS-11-06854.}
\affiliation{Beuth University of Applied Sciences Berlin and\break
University of California, Los Angeles}
\address[a]{Department II---Mathematics, Physics, Chemistry\\
Beuth University of Applied Sciences Berlin\\
Luxemburger Str. 10\\
13353 Berlin\\
Germany\\
\printead{e1}\\
\printead{u1}}

\address[b]{Department of Statistics\\
University of California, Los Angeles\\
8125 Math Sciences Bldg.\\
Box 951554\\
Los Angeles, California 90095-1554\\
USA\\
\printead{e2}\\
\printead{u2}}
\end{aug}

\received{\smonth{7} \syear{2013}}
\revised{\smonth{12} \syear{2013}}

\begin{abstract}
The generalized word length pattern of an orthogonal array allows a ranking
of orthogonal arrays in terms of the generalized minimum aberration
criterion (Xu and Wu
[\textit{Ann. Statist.} \textbf{29} (2001) 1066--1077]). We provide a statistical interpretation for the
number of shortest words of an orthogonal array in terms of sums of $R^{2}$
values (based on orthogonal coding) or sums of squared canonical
correlations (based on arbitrary coding). Directly related to these
results, we derive two versions of generalized resolution for qualitative
factors, both of which are generalizations of the generalized resolution by
Deng and Tang [\textit{Statist. Sinica} \textbf{9} (1999) 1071--1082]
and Tang and Deng [\textit{Ann. Statist.} \textbf{27} (1999) 1914--1926]. We provide a sufficient
condition for one of these to attain its upper bound, and we provide
explicit upper bounds for two classes of symmetric designs. Factor-wise
generalized resolution values provide useful additional detail.
\end{abstract}


\end{frontmatter}

\section{Introduction}\label{sec1}

Orthogonal arrays (OAs) are widely used for designing experiments. One of
the most important criteria for assessing the usefulness of an array is the
generalized word length pattern (GWLP) as proposed by \mbox{\citet{XuWu01}}:
$A_{3}, A_{4}, \ldots$ are the numbers of (generalized) words of lengths
$3, 4, \ldots,$ and the design has resolution $R$, if $A_{i} = 0$ for all $i
< R$ and $A_{R} > 0$. Analogously to the well-known minimum aberration
criterion for regular fractional factorial designs [\citet{FriHun80}],
the quality criterion based on the GWLP is generalized minimum aberration
[GMA; \citet{XuWu01}]: a design $D_{1}$ has better generalized aberration
than a design $D_{2}$, if its resolution is higher or---if both designs
have resolution $R$---if its number $A_{R}$ of shortest words is smaller;
in case of ties in $A_{R}$, frequencies of successively longer words are
compared, until a difference is encountered.

The definition of the $A_{i}$ in Xu and Wu is very technical (see Section~\ref{sec2}). One of the key results of this paper is to provide a statistical
meaning for the number of shortest words, $A_{R}$: we will show that
$A_{R}$ is the sum of $R^{2}$ values from linear models with main effects
model matrix columns in orthogonal coding as dependent variables and full
models in $R -1$ other factors on the explanatory side. For arbitrary
factor coding, the ``sum of $R^{2}$'' interpretation cannot be upheld, but
it can be shown that $A_{R}$ is the sum of squared canonical correlations
[\citet{Hot36}] between a factor's main effects model matrix columns in
arbitrary coding and the full model matrix from $R -1$ other factors. These
results will be derived in Section~\ref{sec2}.

For regular fractional factorial 2-level designs, the GWLP coincides with
the well-known word length pattern (WLP). An important difference between
regular and nonregular designs is that factorial effects in regular
fractional factorial designs are either completely aliased or not aliased
at all, while nonregular designs can have partial aliasing, which can lead
to noninteger entries in the GWLP. In fact, the absence of complete
aliasing has been considered an advantage of nonregular designs [e.g.,
those by Plackett and Burman (\citeyear{PlaBur46})] for screening applications. \citet{DenTan99} and \citet{TanDen99} defined ``generalized resolution''
($\mathrm{GR}$) for nonregular designs with 2-level factors, in order to capture
their advantage over complete confounding in a number. For example, the 12
run Plackett--Burman design has $\mathrm{GR}=3.67$, which indicates that it is
resolution III, but does not have any triples of factors with complete
aliasing. \citet{Evaetal05} have made a useful proposal for
generalizing $\mathrm{GR}$ (called \textit{GRes} by them) for designs in
quantitative factors at 3 levels; in conjunction with \citet{CheYe04},
their proposal can easily be generalized to cover designs with quantitative
factors in general. However, there is so far no convincing proposal for
designs with qualitative factors. The second goal of this paper is to close
this gap, that is, to generalize Deng and Tang's/Tang and Deng's $\mathrm{GR}$ to
OAs for qualitative factors. Any reasonable generalization of $\mathrm{GR}$ has to
fulfill the following requirements: (i)~it must be coding-invariant, that
is, must not depend on the coding chosen for the experimental factors (this
is a key difference vs. designs for quantitative factors), (ii) it must be
applicable for symmetric and asymmetric designs (i.e., designs with a fixed
number of levels and designs with mixed numbers of levels), (iii) like in
the 2-level case, $R+1 > \mathrm{GR} \geq R$ must hold, and $\mathrm{GR} = R$ must be
equivalent to the presence of complete aliasing somewhere in the design,
implying that $R+1 > \mathrm{GR} > R$ indicates a resolution $R$ design with no
complete aliasing among projections of $R$ factors. We offer two proposals
that fulfill all these requirements and provide a rationale behind each of
them, based on the relation of the GWLP to regression relations and
canonical correlations among the columns of the model matrix.

The paper is organized as follows: Section~\ref{sec2} formally introduces
the GWLP and provides a statistical meaning to its number of shortest
words, as discussed above. Section~\ref{sec3} briefly introduces generalized
resolution by Deng and Tang (\citeyear{DenTan99})
and Tang and Deng (\citeyear{DenTan99})
and generalizes it in two meaningful ways. Section~\ref{sec4} shows weak
strength $R$ [in a version modified from \citet{Xu03} to imply strength
$R- 1$] to be sufficient for maximizing one of the generalized
resolutions in a resolution $R$ design. Furthermore, it derives an explicit
upper bound for the proposed generalized resolutions for two classes of
symmetric designs. Section~\ref{sec5} derives factor wise versions of both types of
generalized resolution and demonstrates that these provide useful
additional detail to the overall values. The paper closes with a discussion
and an outlook on future work.

Throughout the paper, we will use the following notation: An orthogonal
array of resolution $R = \mathrm{strength}\ R -1$ in $N$ runs with $n$ factors will
be denoted as $\operatorname{OA}(N, s_{1}, \ldots, s_{n}, R -1)$, with $s_{1}, \ldots,
s_{n}$ the numbers of levels of the $n$ factors (possibly but not
necessarily distinct), or as $\operatorname{OA}(N, s_{1}^{n_1},\ldots, s_{k}^{n_k}$, $R
-1$) with $n_{1}$ factors at $s_{1}$ levels, $\ldots, n_{k}$ factors at
$s_{k}$ levels ($s_{1}, \ldots, s_{k}$ possibly but not necessarily
distinct), whichever is more suitable for the purpose at hand. A subset of
$k$ indices that identifies a $k$-factor projection is denoted by
$\{u_{1},\ldots,u_{k} \} (\subseteq \{1,\ldots,n\})$. The unsquared letter
$R$ always refers to the resolution of a design, while $R^{2}$ denotes the
coefficient of determination.

\section{Projection frequencies and linear models}\label{sec2}

Consider an $\operatorname{OA}(N, s_{1}, \ldots,\break  s_{n}, R -1)$. The resolution $R$
implies that main effects can be confounded with interactions among $R -1$
factors, where the extent of confounding of degree $R$ can be investigated
on a global scale or in more detail: Following \citet{XuWu01}, the
factors are coded in orthogonal contrasts with squared column length
normalized to $N$. We will use the expression ``normalized orthogonal
coding'' to refer to this coding; on the contrary, the expressions
``orthogonal coding'' or ``orthogonal contrast coding'' refer to main
effects model matrix columns that have mean zero and are pairwise
orthogonal, but need not be normalized. For later reference, note that for
orthogonal coding (whether normalized or not) the main effects model matrix
columns for an OA (of strength at least 2) are always uncorrelated.

We write the model matrix for the full model in normalized orthogonal
coding as
%
\begin{equation}
\mathbf{M} = (\mathbf{M}_{0}, \mathbf{M}_{1}, \ldots,
\mathbf{M}_{n}), \label{eq1}
\end{equation}
where $\mathbf{M}_{0}$ is a column of ``$+1$''s, $\mathbf{M}_{1}$ contains
all main effects model matrices, and $\mathbf{M}_{k}$ is the matrix of all
$ { n \choose
k}$ $k$-factor interaction model matrices, $k = 2,\ldots,n$.
The portion $\mathbf{X}_{u_1,\ldots, u_k}$ of $\mathbf{M}_{k} =
(\mathbf{X}_{1,\ldots, k},\ldots,\mathbf{X}_{n - k+1,\ldots, n})$ denotes the
model matrix for the particular $k$-factor interaction indexed by
$\{u_{1},\ldots,u_{k}\}$ and is obtained by all products from one main
effects contrast column each from the $k$ factors in the interaction. Note
that the normalized orthogonal coding of the main effects implies that all
columns of $\mathbf{M}_{k}$ have squared length $N$ for $k \leq R -1$. Now,
on the global scale, the overall number of words of length $k$ can be
obtained as the sum of squared column averages of $\mathbf{M}_{k}$, that
is, $A_{k} = \mathbf{1}_{N}^{\mathrm{T}} \mathbf{M}_{k}
\mathbf{M}_{k}^{\mathrm{T}} \mathbf{1}_{N} / N^{2}$. Obviously, this sum
can be split into contributions from individual $k$-factor projections for
more detailed considerations, that is,
%
\begin{equation}
A_{k} = \mathop{\sum_{\{ u_{1},\ldots,u_{k}  \} }}_{\subseteq  \{ 1,\ldots,n  \}} \mathbf{1}_{N}^{\mathrm{T}}
\mathbf{X}_{u_{1},\ldots,u_{k}}\mathbf{X}_{u_{1},\ldots, u_{k}}^{\mathrm{T}}
\mathbf{1}_{N} / N^{2} =:\mathop{\sum_{ \{ u_{1},\ldots,u_{k}  \}}}_{
 \subseteq  \{ 1,\ldots,n  \}}
a_{k} ( u_{1},\ldots,u_{k} ),\label{eq2}
\end{equation}
where $a_{k}(u_{1},\ldots,u_{k})$ is simply the $A_{k}$ value of the
$k$-factor projection $\{u_{1},\ldots,\break u_{k}\}$. The summands
$a_{k}(u_{1},\ldots,u_{k})$ are called ``projection frequencies.''

\begin{example}\label{ex1}
For 3-level factors, normalized polynomial coding
has the linear contrast coefficients $- \sqrt{3 / 2}, 0, \sqrt{3 / 2}$ and
the quadratic contrast coefficients $\sqrt{1 / 2}, - \sqrt{2}, \sqrt{1
/ 2}$. For the regular design $\operatorname{OA}(9, 3^{3}, 2)$ with the defining relation
$\mathrm{C}=\mathrm{A}+\mathrm{B}$ (mod 3), the model matrix $\mathbf{M}$
has dimensions $9\times27$, including one column for $\mathbf{M}_{0}$, six for
$\mathbf{M}_{1}$, twelve for $\mathbf{M}_{2}$ and eight for
$\mathbf{M}_{3}$. Like always, the column sum of $\mathbf{M}_{0}$ is $N$
(here: 9), and like for any orthogonal array, the column sums of
$\mathbf{M}_{1}$ and $\mathbf{M}_{2}$ are 0, which implies $A_{0} =1$,
$A_{1} =A_{2} =0$. We now take a closer look at $\mathbf{M}_{3}$, arranging
factor A as (0 0 0 1 1 1 2 2 2), factor B as (0 1 2 0 1 2 0 1 2) and factor
C as their sum (mod 3), denoting linear contrast columns by the subscript
$l$ and quadratic contrast columns by the subscript $q$. Then

{\fontsize{8}{10}\selectfont{
\[
\hspace*{-6pt}\mathbf{M}_{3} = \pmatrix{ \mathrm{contrast} & \mathrm{A}_{
l}
\mathrm{B}_{ l}\mathrm{C}_{ l} & \mathrm{A}_{ q}
\mathrm{B}_{
l}\mathrm{C}_{ l} & \mathrm{A}_{ l}
\mathrm{B}_{ q}\mathrm{C}_{ l} & \mathrm{A}_{ q}
\mathrm{B}_{ q}\mathrm{C}_{ l} & \mathrm{A}_{ l}
\mathrm{B}_{
l}\mathrm{C}_{ q} & \mathrm{A}_{ q}
\mathrm{B}_{ l}\mathrm{C}_{ q} & \mathrm{A}_{ l}
\mathrm{B}_{ q}\mathrm{C}_{ q} & \mathrm{A}_{ q}
\mathrm{B}_{
q}\mathrm{C}_{ q} \vspace*{2pt}
\cr
\hline & -
\sqrt{\frac{27}{8}} & \sqrt{\frac{9}{8}} & \sqrt{\frac{9}{8}} & -
\sqrt{\frac{3}{8}} & \sqrt{\frac{9}{8}} & - \sqrt{\frac{3}{8}}
& - \sqrt{\frac{3}{8}} & \sqrt{\frac{1}{8}} \vspace*{2pt}
\cr
& 0 & 0
& 0 & 0 & 0 & 0 & - \sqrt{6} & \sqrt{2} \vspace*{2pt}
\cr
& - \sqrt{
\frac{27}{8}} & \sqrt{\frac{9}{8}} & - \sqrt{\frac{9}{8}} &
\sqrt{\frac{3}{8}} & - \sqrt{\frac{9}{8}} & \sqrt{\frac{3}{8}} &
- \sqrt{\frac{3}{8}} & \sqrt{\frac{1}{8}} \vspace*{2pt}
\cr
& 0 & 0 &
0 & 0 & 0 & - \sqrt{6} & 0 & \sqrt{2} \vspace*{2pt}
\cr
& 0 & 0 & 0 & \sqrt{6} & 0
& 0 & 0 & \sqrt{2} \vspace*{2pt}
\cr
& 0 & \sqrt{\frac{9}{2}} & 0 & \sqrt{
\frac{3}{2}} & 0 & - \sqrt{\frac{3}{2}} & 0 & - \sqrt{
\frac{1}{2}} \vspace*{2pt}
\cr
& - \sqrt{\frac{27}{8}} & - \sqrt{
\frac{9}{8}} & \sqrt{\frac{9}{8}} & \sqrt{\frac{3}{8}} & -
\sqrt{\frac{9}{8}} & - \sqrt{\frac{3}{8}} & \sqrt{\frac{3}{8}} &
\sqrt{\frac{1}{8}} \vspace*{2pt}
\cr
& 0 & 0 & \sqrt{\frac{9}{2}} &
\sqrt{\frac{3}{2}} & 0 & 0 & - \sqrt{\frac{3}{2}} & - \sqrt{
\frac{1}{2}} \vspace*{2pt}
\cr
& 0 & 0 & 0 & 0 & - \sqrt{\frac{9}{2}} &
- \sqrt{\frac{3}{2}} & - \sqrt{\frac{3}{2}} & - \sqrt{
\frac{1}{2}} \vspace*{2pt}
\cr
\hline
\mathrm{column}  & -
\sqrt{\frac{243}{8}} & \sqrt{\frac{81}{8}} & \sqrt{\frac{81}{8}} &
\sqrt{\frac{243}{8}} & - \sqrt{\frac{81}{8}} & - \sqrt{\frac{243}{8}}
& - \sqrt{\frac{243}{8}} & \sqrt{\frac{81}{8}}
\vspace*{2pt}\cr
\mathrm{sum}&&&&&&&& }.
\]}}
\hspace*{-3pt}Half of the squared column sums of $\mathbf{M}_{3}$ are $243/8$ and $81/8$,
respectively. This implies that the sum of the squared column sums is
$A_{3} = a_{3} (1,2,3) = (4\cdot 243/8+4\cdot 81/8)/81 = 2$.
\end{example}

\begin{example}\label{ex2}
Table~\ref{tab1} displays the only $\operatorname{OA}(18, 2^{1}3^{2}, 2)$
that cannot be obtained as a projection from the L18 design that was
popularized by Taguchi (of course, this triple is not interesting as a
stand-alone design, but as a projection from a design in more factors only;
the Taguchi L18 is for convenience displayed in Table~\ref{tab4} below). The 3-level
factors are coded like in Example~\ref{eq1}, for the 2-level factor, normalized
orthogonal coding is the customary $-1/+1$ coding. Now, the model matrix
$\mathbf{M}$ has dimensions $18\times18$, including one column for
$\mathbf{M}_{0}$, five for $\mathbf{M}_{1}$, eight for $\mathbf{M}_{2}$ and
four for $\mathbf{M}_{3}$. Again, $A_{0} =1$, $A_{1} =A_{2} =0$. The
squared column sums of $\mathbf{M}_{3}$ are 9 ($1\times$), 27 ($2\times$) and 81 ($1\times$),
respectively. Thus, $A_{3} = a_{3} (1,2,3) = (9+2\cdot 27+81)/324 = 4/9$.
\end{example}

The projection frequencies $a_{k}(u_{1},\ldots,u_{k})$ from equation (\ref{eq2}) are
the building blocks for the overall $A_{k}$. The $a_{R}(u_{1},\ldots,u_{R})$
will be instrumental in defining one version of generalized resolution.
Theorem~\ref{th1} provides them with an intuitive interpretation. The proof is
given in the \hyperref[app]{Appendix}.

\begin{table}
\caption{A partially confounded $\operatorname{OA}(18, 2^{1}3^{2}, 2)$
(transposed)}\label{tab1}
\begin{tabular*}{305pt}{@{\extracolsep{\fill}}lcccccccccccccccccc@{}}
\hline
A &0 &1 &1 &0& 0 &1 &0& 1 &0& 1& 1& 0& 1& 0& 0& 1& 0& 1\\
B& 0& 0 &0& 0 &0& 0& 1 &1 &1 &1 &1& 1 &2 &2& 2& 2 &2 &2\\
C& 0& 1& 2& 1& 2& 0 &0 &2& 0& 1& 1& 2& 2& 1& 2& 0&1 &0\\
\hline
\end{tabular*}  \vspace*{-3pt}
\end{table}

\begin{theorem}\label{th1}
In an $\operatorname{OA}(N, s_{1}, \ldots, s_{n}, R -1)$, denote
by $\mathbf{X}_{c}$ the model matrix for the main effects of a particular
factor $c \in \{u_{1},\ldots,u_{R} \} \subseteq \{1,\ldots,n\}$ in normalized
orthogonal coding, and let $\mathrm{C} = \{u_{1},\ldots,u_{R}\} \setminus
\{c\}$. Then $a_{R}(u_{1},\ldots,\break u_{R})$ is the sum of the $R^{2}$-values
from the $s_{c} -1$ regression models that explain the columns of
$\mathbf{X}_{c}$ by a full model in the factors from $\mathrm{C}$.
\end{theorem}

\begin{remark}\label{re1}
(i) Theorem~\ref{th1} holds regardless which factor is
singled out for the left-hand side of the model. (ii) The proof simplifies
by restriction to normalized orthogonal coding, but the result holds
whenever the factor $c$ is coded by any set of orthogonal contrasts,
whether normalized or not. (iii) Individual $R^{2}$ values are coding
dependent, but the sum is not. (iv) In case of normalized orthogonal coding
for all factors, the full model in the factors from $\mathrm{C}$ can be reduced to the
$R -1$ factor interaction only, since the matrix $\mathbf{X}_{c}$ is
orthogonal to the model matrices for all lower degree effects in the other
$R -1$ factors.
\end{remark}

\setcounter{example}{0}
\begin{example}[(Continued)]\label{ex1co} The overall $a_{3} (1,2,3) = 2$ is the sum
of two $R^{2}$ values which are 1, regardless\vadjust{\goodbreak} which factor is singled out
as the main effects factor for the left-hand sides of regression. This
reflects that the level of each factor is uniquely determined by the level
combination of the other two factors.
\end{example}

\begin{example}[(Continued)]\label{ex2co} The $R^{2}$ from regressing the single
model matrix column of the 2-level factor on the four model matrix columns
for the interaction among the two 3-level factors is $4/9$. Alternatively,
the $R^{2}$-values for the regression of the two main effects columns for
factor B on the AC interaction columns are $1/9$ and $3/9$, respectively, which
also yields the sum $4/9$ obtained above for $a_{3}(1,2,3)$. For factor B in
dummy coding with reference level 0 instead of normalized polynomical
coding, the two main effects model matrix columns for factor B have
correlation 0.5; the sum of the $R^{2}$ values from full models in A and C
for explaining these two columns is $1/3 + 1/3 = 2/3 \neq a_{3} (1,2,3) =
4/9$. This demonstrates that Theorem~\ref{th1} is not applicable if orthogonal
coding [see Remark~\ref{re1}(ii)] is violated.
\end{example}

\begin{corollary}\label{co1} In an $\operatorname{OA}(N, s_{1}, \ldots,s_{n}, R -1)$, let
$\{u_{1},\ldots,u_{R} \} \subseteq \{1,\ldots,\break  n\}$, with $s_{\mathrm{min}}=
\mathrm{min}_{i=1,\ldots,R}(s_{u_i})$.

\begin{longlist}[(iii)]
\item[(i)] A factor $c \in \{u_{1},\ldots,u_{R}\}$ in $s_{c}$ levels is
completely confounded by the factors in
$\mathrm{C} = \{u_{1},\ldots,u_{R}\} \setminus \{ c\}$, if and only if
$a_{R}(u_{1},\ldots,u_{R}) = s_{c} -1$.
\item[(ii)] $a_{R}(u_{1},\ldots,u_{R}) \leq s_{\mathrm{min}} - 1$.
\item[(iii)] If several factors in $\{u_{1},\ldots,u_{R}\}$ have
$s_{\mathrm{min}}$ levels, either all of them are or none of them is completely
confounded by the respective other $R -1$ factors in
$\{u_{1},\ldots,u_{R}\}$.
\item[(iv)] A factor with more than $s_{\mathrm{min}}$ levels cannot be
completely confounded by the other factors in $\{u_{1},\ldots,u_{R}\}$.
\end{longlist}
\end{corollary}

Part (i) of Corollary~\ref{co1} follows easily from Theorem~\ref{th1}, as
$a_{R}(u_{1},\ldots,u_{R}) = s_{c} -1$ if and only if all $R^{2}$ values for
columns of the factor $c$ main effects model matrix are 100\%, that is, the
factor $c$ main effects model matrix is completely explained by the factors
in C. Part (ii) follows, because the sum of $R^{2}$ values is of course
bounded by the minimum number of regressions conducted for any single
factor $c$, which is $s_{\mathrm{min}}- 1$. Parts (iii) and (iv) follow
directly from parts (i) and (ii). For symmetric $s$-level designs, part
(ii) of the corollary has already been proven by \citet{XuCheWu04}.

\begin{table}
\caption{An $\operatorname{OA}(8, 4^{1}2^{2}, 2)$ (transposed)}\label{tab2}
\begin{tabular*}{150pt}{@{\extracolsep{\fill}}lcccccccc@{}}
\hline
A &0 &0 &0& 0 &1 &1 &1 &1
\\
B &0 &0& 1& 1 &0& 0 &1 &1
\\
C &0 &2& 1& 3& 3& 1& 2 &0
\\
\hline
\end{tabular*}
\end{table}

\begin{example}\label{ex3}
For the design of Table~\ref{tab2}, $s_{\mathrm{min}} = 2$, and
$a_{3} (1,2,3) = 1$, that is, both 2-level factors are completely
confounded, while the 4-level factor is only partially confounded. The
individual $R^{2}$ values for the separate degrees of freedom of the
4-level factor main effect model matrix depend on the coding (e.g., 0.2, 0
and 0.8 for the linear, quadratic and cubic contrasts in normalized
orthogonal polynomial coding), while their sum is 1, regardless of the
chosen orthogonal coding.
\end{example}

\begin{theorem}\label{th2}
In an $\operatorname{OA}(N, s_{1}, \ldots, s_{n}, R -1)$, let
$\{u_{1},\ldots,u_{R} \} \subseteq \{1,\ldots,n\}$ with $s_{\mathrm{min}} =
\mathrm{min}_{i=1,\ldots,R}(s_{u_i})$. Let $c \in \{u_{1},\ldots,u_{R}\}$ with $s_{c}
= s_{\mathrm{min}}$, $\mathrm{C} = \{u_{1},\ldots,\break  u_{R}\} \setminus \{c\}$.
Under normalized orthogonal coding denote by $\mathbf{X}_{c}$ the main
effects model matrix for factor $c$ and by $\mathbf{X}_{\mathrm{C}}$ the $R
-1$ factor interaction model matrix for the factors in $\mathrm{C}$.

If $a_{R}(u_{1},\ldots,u_{R}) = s_{\mathrm{min}} - 1$, $\mathbf{X}_{\mathrm{C}}$
can be orthogonally transformed (rotation and or switching) such that
$s_{\mathrm{min}} -1$ of its columns are collinear to the columns of~$\mathbf{X}_{c}$.
\end{theorem}

\begin{pf}
$a_{R}(u_{1},\ldots,u_{R}) = s_{\mathrm{min}} - 1$ implies
all $s_{\mathrm{min}}- 1$ regressions of the columns of $\mathbf{X}_{c}$ on
the columns of $\mathbf{X}_{\mathrm{C}}$ have $R^{2} =1$. Then, each of the
$s_{\mathrm{min}} - 1$ $\mathbf{X}_{c}$ columns can be perfectly matched by a
linear combination $\mathbf{X}_{\mathrm{C}} \mathbf{b}$ of the
$\mathbf{X}_{\mathrm{C}}$ columns; since all columns have the same length,
this linear transformation involves rotation and/or switching only. If
necessary, these $s_{\mathrm{min}} - 1$ orthogonal linear combinations can be
supplemented by further length-preserving orthogonal linear combinations so
that the dimension of $\mathbf{X}_{\mathrm{C}}$ remains intact.
\end{pf}

Theorems \ref{th1} and \ref{th2} are related to canonical correlation analysis, and the
redundancy index discussed in that context [\citet{SteLov68}]. In
order to make the following comments digestible, a brief definition of
canonical correlation analysis is included without going into any technical
detail about the method; details can, for example, be found in H\"{a}rdle and
Simar [(\citeyear{HarSim03}), Chapter~14]. It will be helpful to think of the columns of the
main effects model matrix of factor $c$ as the $Y$ variables and the
columns of the full model matrix in the $R -1$ other factors from the set $\mathrm{C}$
(excluding the constant column of ones for the intercept) as the $X$
variables of the following definition and explanation. As it would be
unnatural to consider the model matrices from experimental designs as
random variables, we directly define canonical correlation analysis in
terms of data matrices $\mathbf{X}$ and $\mathbf{Y}$ ($N$ rows each) and
empirical covariance matrices $\mathbf{S}_{xx} =
\mathbf{X}^{*\mathrm{T}} \mathbf{X}^{*} /(N - 1)$, $\mathbf{S}_{yy} =
\mathbf{Y}^{*\mathrm{T}} \mathbf{Y}^{*} /(N - 1)$, $\mathbf{S}_{xy}=
\mathbf{X}^{*\mathrm{T}} \mathbf{Y}^{*} /(N - 1)$ and $\mathbf{S}_{yx}
= \mathbf{Y}^{*\mathrm{T}} \mathbf{X}^{*} /(N - 1)$, where the superscript
$^{\ast}$ denotes columnwise centering of a matrix. We do not attempt a minimal
definition, but prioritize suitability for our purpose. Note that our
$\mathbf{S}_{xx}$ and $\mathbf{S}_{yy}$ are nonsingular matrices,
since the designs we consider have strength $R -1$; the covariance matrix
$(\mathbf{X}^{*} \mathbf{Y}^{*})^{\mathrm{T}}(\mathbf{X}^{*}
\mathbf{Y}^{*})/(N - 1)$ of the combined set of variables may, however, be
singular, which does not pose a problem to canonical correlation analysis,
even though some accounts request this matrix to be nonsingular.

\begin{definition}\label{de1}
Consider a set of $p X$-variables and $q
Y$-variables. Let the $N\times p$ matrix $\mathbf{X}$ and the
$N\times q$ matrix $\mathbf{Y}$ denote the data matrices of $N$
observations, and $\mathbf{S}_{xx}$, $\mathbf{S}_{yy}$,
$\mathbf{S}_{xy}$ and $\mathbf{S}_{yx}$ the empirical covariance
matrices obtained from them, with positive definite $\mathbf{S}_{xx}$
and $\mathbf{S}_{yy}$.
\begin{enumerate}[(ii)]
\item[(i)] Canonical correlation analysis creates $k = \mathrm{min}(p, q)$ pairs
of linear combination vectors $\mathbf{u}_{i} =\mathbf{Xa}_{i}$ and
$\mathbf{v}_{i} =\mathbf{Yb}_{i}$ with $p\times 1$ coefficient vectors
$\mathbf{a}_{i}$ and $q\times 1$ coefficient vectors $\mathbf{b}_{i}$, $i =
1,\ldots,k$, such that:
\begin{enumerate}[(a)]
\item[(a)] the $\mathbf{u}_{1},\ldots,\mathbf{u}_{k}$ are uncorrelated to
each other,
\item[(b)] the $\mathbf{v}_{1},\ldots,\mathbf{v}_{k}$ are uncorrelated to
each other,\vadjust{\goodbreak}
\item[(c)] the pair ($\mathbf{u}_{1}$, $\mathbf{v}_{1}$) has the maximum
possible correlation for any pair of linear combinations of the
$\mathbf{X}$ and $\mathbf{Y}$ columns, respectively,
\item[(d)] the pairs ($\mathbf{u}_{i}$, $\mathbf{v}_{i}$), $i=2,\ldots,k$
successively maximize the remaining correlation, given the constraints of
(a) and (b).
\end{enumerate}
\item[(ii)] The correlations $r_{i} = \operatorname{cor}(\mathbf{u}_{i},
\mathbf{v}_{i})$ are called ``canonical correlations,''
and the $\mathbf{u}_{i}$ and $\mathbf{v}_{i}$ are called ``canonical
variates.''
\end{enumerate}
\end{definition}

\begin{remark}\label{re2}
(i) If the matrices $\mathbf{X}$
and $\mathbf{Y}$ are centered, that is,
$\mathbf{X}=\mathbf{X}^{*}$ and
$\mathbf{Y}=\mathbf{Y}^{*}$, the
$\mathbf{u}$ and $\mathbf{v}$
vectors also have zero means, and the uncorrelatedness in~(a) and
(b) is equivalent to orthogonality of the vectors. (ii) It is well known
that the canonical correlations are the eigenvalues of the matrices
$\mathbf{Q}_{1}=\mathbf{S}_{xx}^{-1}\mathbf{S}_{xy}\mathbf{S}_{yy}^{-1}\mathbf{S}_{yx}$
and $\mathbf{Q}_{2}=\mathbf{S}_{yy}^{-1}\mathbf{S}_{yx}\mathbf{S}_{xx}^{-1}\mathbf{S}_{xy}$
[the first $\mathrm{min}({p},
{q})$ eigenvalues of both matrices are the
same; the larger matrix has the appropriate number of additional zeroes]
and the $\mathbf{a}_{i}$
are the corresponding eigenvectors of
$\mathbf{Q}_{1}$, the
$\mathbf{b}_{i}$ the
corresponding eigenvectors of
$\mathbf{Q}_{2}$.
\end{remark}

According to the definition, the canonical correlations are nonnegative.
It can also be shown that $\mathbf{u}_{i}$ and $\mathbf{v}_{j}$, $i \neq
j$, are uncorrelated, and orthogonal in case of centered data matrices;
thus, the pairs ($\mathbf{u}_{i}$, $\mathbf{v}_{i}$) decompose the relation
between $\mathbf{X}$ and $\mathbf{Y}$ into uncorrelated components, much
like the principal components decompose the total variance into
uncorrelated components. In data analysis, canonical correlation analysis
is often used for dimension reduction. Here, we retain the full
dimensionality. For uncorrelated $Y$ variables like the model matrix
columns of $\mathbf{X}_{c}$ in Theorem~\ref{th1}, it is straightforward to see that
the sum of the $R^{2}$ values from regressing each of the $Y$ variables on
all the $X$ variables coincides with the sum of the squared canonical
correlations. It is well known that the canonical correlations are
invariant to arbitrary nonsingular affine transformations applied to the
$X$- and $Y$-variables, which translate into nonsingular linear
transformations applied to the centered $\mathbf{X}$- and
$\mathbf{Y}$-matrices [cf., e.g., H\"{a}rdle and Simar (\citeyear{HarSim03}),
Theorem~14.3].
For our application, this implies invariance of the canonical correlations
to factor coding. Unfortunately, this invariance property does not hold for
the $R^{2}$ values or their sum: according to Lazraq and Cl\'{e}roux
[(\citeyear{LazCle01}), Section~2] the aforementioned redundancy index---which is the average
$R^{2}$ value calculated as $a_{R}(u_{1},\ldots,u_{R})/(s_{c} -1)$ in the
situation of Theorem~\ref{th1}---is invariant to linear transformations of the
centered $\mathbf{X}$ matrix, but only to \textit{orthonormal}
transformations of the centered $\mathbf{Y}$ matrix or scalar multiples
thereof. For correlated $Y$-variables, the redundancy index contains some
overlap between variables, as was already seen for Example \ref{ex2}, where the sum
of the $R^{2}$ values from dummy coding exceeded $a_{3}(1,2,3)$; in that
case, only the average or sum of the squared canonical correlations yields
an adequate measure of the overall explanatory power of the $X$-variables
on the $Y$-variables. Hence, for the case of arbitrary coding, Theorem~\ref{th1}
has to be restated in terms of squared canonical correlations.

\begin{theorem}\label{th3}
In an $\operatorname{OA}(N, s_{1},\ldots, s_{n}, R -1)$, denote
by $\mathbf{X}_{c}$ the model matrix for the main effects of a particular
factor\vadjust{\goodbreak} $c \in \{u_{1},\ldots,u_{R}\}$ in arbitrary coding, and let
$\mathrm{C} = \{u_{1},\ldots,u_{R}\}\setminus \{c\}$. Then
$a_{R}(u_{1},\ldots,u_{R})$ is the sum of the squared canonical correlations
from a canonical correlation analysis of the columns of $\mathbf{X}_{c}$
and the columns of the full model matrix $\mathbf{F}_{\mathrm{C}}$ in the
factors from $\mathrm{C}$.
\end{theorem}

\setcounter{example}{0}
\begin{example}[(Continued)]
$s_{\mathrm{min}}=3$, $a_{3} (1,2,3) = 2$,
that is, the assumptions of Theorems \ref{th2} and \ref{th3} are fulfilled. Both canonical
correlations must be 1, because the sum must be 2. The transformation of
$\mathbf{X}_{\mathrm{C}}$ from Theorem~\ref{th2} can be obtained from the canonical
correlation analysis: For all factors in the role of $Y$, $\mathbf{v}_{i}
\propto \mathbf{y}_{i}$ (with $\mathbf{y}_{i}$ denoting the $i$th column
of the main effects model matrix of the $Y$-variables factor) can be used.
For the first or second factor in the role of $Y$, the corresponding
canonical vectors on the $X$ side fulfill
\begin{eqnarray}
\mathbf{u}_{1}& \propto& \mathrm{B}_{q}\mathrm{C}_{l} -
\mathrm{B}_{l}\mathrm{C}_{q} - \sqrt{3}\mathrm{B}_{l}\mathrm{C}_{l} -
\sqrt{3}\mathrm{B}_{q}\mathrm{C}_{q},\nonumber\\
\mathbf{u}_{2} &\propto& \sqrt{3}\mathrm{B}_{l}\mathrm{C}_{q} -
\sqrt{3}\mathrm{B}_{q}\mathrm{C}_{l} - \mathrm{B}_{l}\mathrm{C}_{l} -
\mathrm{B}_{q}\mathrm{C}_{q}\nonumber\\
\eqntext{\mbox{(or B replaced by A for the second factor in
the role of $Y$),}}
\end{eqnarray}
with the indices $l$ and $q$ denoting the
\textit{normalized} linear and quadratic coding introduced above.
For the third factor in the role of $Y$,
\begin{eqnarray*}
\mathbf{u}_{1}& \propto& - \sqrt{3}\mathrm{A}_{l}\mathrm{B}_{l} +
\mathrm{A}_{q}\mathrm{B}_{l} + \mathrm{A}_{l}\mathrm{B}_{q} +
\sqrt{3}\mathrm{A}_{q}\mathrm{B}_{q},\\
\mathbf{u}_{2}& \propto& -A_{l}\mathrm{B}_{l} -
\sqrt{3}\mathrm{A}_{l}\mathrm{B}_{q} - \sqrt{3}\mathrm{A}_{q}\mathrm{B}_{l}
+ \mathrm{A}_{q}\mathrm{B}_{q}.
\end{eqnarray*}
\end{example}

\setcounter{example}{0}
\begin{example}[(Now with dummy coding)]
When using the design of
Example~\ref{ex1} for an experiment with qualitative factors, dummy coding is much
more usual than orthogonal contrast\vadjust{\goodbreak} coding. This example shows how Theorem~\ref{th3} can be applied for arbitrary nonorthogonal coding: $\mathrm{A}_{1}$ is 1
for $\mathrm{A}=1$ and 0 otherwise, $\mathrm{A}_{2}$ is 1 for
$\mathrm{A}=2$ and 0 otherwise, B and C are coded analogously; interaction
matrix columns are obtained as products of the respective main effects
columns. The main effect and two-factor interaction model matrix columns in
this coding do not have column means zero and have to be centered first by
subtracting $1/3$ or $1/9$, respectively. As canonical correlations are
invariant to affine transformations, dummy coding leads to the same
canonical correlations as the previous normalized orthogonal polynomial
coding. We consider the first factor in the role of $Y$; the centered model
matrix columns $\mathbf{y}_{1} = \mathrm{A}_{1} -1/3$ and $\mathbf{y}_{2} =
\mathrm{A}_{2} -1/3$ are correlated, so that we must not choose both
canonical variates for the $Y$ side proportional\vspace*{1pt} to the original variates.
One instance of the canonical variates for the $Y$ side is $\mathbf{v}_{1}
= - \mathbf{y}_{1} / \sqrt{2}, \mathbf{v}_{2} =  ( \mathbf{y}_{1} +
2\mathbf{y}_{2}  ) / \sqrt{6}$; these canonical vectors are unique up
to rotation only, because the two canonical correlations have the same
size. The corresponding canonical vectors on the $X$ side are obtained from
the centered full model matrix
\begin{eqnarray*}
&&\mathbf{F}_{\mathrm{C}} = \bigl( \bigl( \mathrm{B}_{1} -
\tfrac{1}{3} \bigr), \bigl( \mathrm{B}_{2} - \tfrac{1}{3}
\bigr), \bigl( \mathrm{C}_{1} - \tfrac{1}{3} \bigr), \bigl(
\mathrm{C}_{2} - \tfrac{1}{3} \bigr), \bigl(
\mathrm{B}_{1}\mathrm{C}_{1} - \tfrac{1}{9} \bigr),
\bigl( \mathrm{B}_{2}\mathrm{C}_{1} - \tfrac{1}{9}
\bigr), \\
&&\hspace*{212pt}{}\bigl( \mathrm{B}_{1}\mathrm{C}_{2} -
\tfrac{1}{9} \bigr), \bigl( \mathrm{B}_{2}\mathrm{C}_{2}
- \tfrac{1}{9} \bigr) \bigr)\vadjust{\goodbreak}
\end{eqnarray*}
as $\mathbf{u}_{1} =  ( - \mathbf{f}_{2} - \mathbf{f}_{3} +
\mathbf{f}_{5} + 2\mathbf{f}_{6} - \mathbf{f}_{7} + \mathbf{f}_{8}  )
/ \sqrt{2}$ and $\mathbf{u}_{2} =  ( 2\mathbf{f}_{1} + \mathbf{f}_{2} +
\mathbf{f}_{3} + 2\mathbf{f}_{4} - 3\mathbf{f}_{5} - 3\mathbf{f}_{7} -
3\mathbf{f}_{8}  ) / \sqrt{6}$, with $\mathbf{f}_{j}$ denoting the
$j$th column of $\mathbf{F}_{\mathrm{C}}$.

Note that the canonical vectors $\mathbf{u}_{1}$ and $\mathbf{u}_{2}$ now
contain contributions not only from the interaction part of the model
matrix but also from the main effects part, that is, we do indeed need the
full model matrix as stated in Theorem~\ref{th3}.
\end{example}

\begin{example}[(Continued)]
$s_{\mathrm{min}} = 2$, $a_{3} (1,2,3) =
4/9$, that is, the assumption of Theorem~\ref{th2} is not fulfilled, the assumption
of Theorem~\ref{th3} is. The canonical correlation using the one column main
effects model matrix of the 2-level factor A in the role of $Y$ is $2/3$, the
canonical correlations using the main effects model matrix for the 3-level
factor B in the role of $Y$ are $2/3$ and 0; in both cases, the sum of the
squared canonical correlations is $a_{3} (1,2,3) = 4/9$. For any other
coding, for example, the dummy coding for factor B considered earlier, the
canonical correlations remain unchanged ($2/3$ and 0, resp.), since
they are coding invariant; thus, the sum of the squared canonical
correlations remains $4/9$, even though the sum of the $R^{2}$ values was
found to be different. Of course, the linear combination coefficients for
obtaining the canonical variates depend on the coding [see, e.g.,
H\"{a}rdle and Simar (\citeyear{HarSim03}), Theorem~14.3].
\end{example}

\begin{table}
\caption{Main effects matrix of factor A regressed on full model in factors
B and C for the 10 nonisomorphic GMA $\operatorname{OA}(32, 4^{3}, 2)$}
\label{tab3}
\begin{tabular*}{\textwidth}{@{\extracolsep{\fill}}lcccccccccc@{}}
\hline
\multicolumn{3}{c}{$\bolds{R^{2}}$ \textbf{values from}} &
\multicolumn{3}{c}{$\bolds{R^{2}}$ \textbf{values from}} &
\multicolumn{3}{c}{\textbf{Squared canonical}} & &
\\
\multicolumn{3}{c}{\textbf{polynomial coding}} &
\multicolumn{3}{c}{\textbf{Helmert coding}} &
\multicolumn{3}{c}{\textbf{correlations}} & &
\\[-6pt]
\multicolumn{3}{@{}l}{\hrulefill} &
\multicolumn{3}{c}{\hrulefill} &
\multicolumn{3}{c}{\hrulefill} & &
\\
\textbf{{L}} & \textbf{{Q}} & \textbf{{C}} &
\textbf{1} & \textbf{2} & \textbf{3} & \textbf{1} & \textbf{2} &
\textbf{3} & $\bolds{A_{3}}$ & \multicolumn{1}{c@{}}{\textbf{Designs}}\\
\hline
0.8 & 0\phantom{00.} & 0.2\phantom{00} & 0 & $2/3$ & $1/3$ & 1\phantom{000.} & 0\phantom{000.} & 0\phantom{00.} & 1 & 1\\
0.65 & 0\phantom{00.} & 0.35\phantom{0} & $1/8$ & $13/24$ & $1/3$ & 0.75\phantom{0} & 0.25\phantom{0} & 0\phantom{00.} & 1 & 2\\
0.5 & 0\phantom{00.} & 0.5\phantom{00} & $1/4$ & $5/12$ & $1/3$ & 0.5\phantom{00} & 0.5\phantom{00} & 0\phantom{00.} & 1 & 3, 6, 8, 10\\
0.45 & 0.25 & 0.3\phantom{00} & $1/4$ & $5/12$ & $1/3$ & 0.5\phantom{00} & 0.25\phantom{0} & 0.25 & 1 & 4, 5, 7\\
0.375 & 0.25 & 0.375 & $5/16$ & $17/48$ & $1/3$ & 0.375 & 0.375 & 0.25 & 1 & 9\\
\hline
\end{tabular*}
\end{table}

Canonical correlation analysis can also be used to verify that a result
analogous to Theorem~\ref{th2} cannot be generalized to sets of $R$ factors for
which $a_{R}(u_{1},\ldots,u_{R}) < s_{\mathrm{min}} - 1$. For this, note that
the number of nonzero\vadjust{\goodbreak} canonical correlations indicates the dimension of
the relationship between the $X$- and the $Y$-variables.

Table~\ref{tab3} displays the $R^{2}$ values from two different orthogonal codings
and the squared canonical correlations from the main effects matrix of the
first factor ($Y$-variables) vs. the full model matrix of the other two
factors \mbox{($X$-variables)} for the ten nonisomorphic GMA $\operatorname{OA}(32,4^{3},2)$
obtained from \citet{EenSch}. These designs have one
generalized word of length 3, that is, they are nonregular. There are cases
with one, two and three nonzero canonical correlations, that is, neither
is it generally possible to collapse the linear dependence into a
one-dimensional structure nor does the linear dependence generally involve
more than one dimension.

\section{Generalized resolution}\label{sec3}

Before presenting the new proposals for generalized resolution, we briefly
review generalized resolution for symmetric 2-level designs by \citet{DenTan99} and \citet{TanDen99}. For 2-level factors, each effect has
a single degree of freedom (df) only, that is, all the $\mathbf{X}$'s in
any $\mathbf{M}_{k}$ [cf. equation (\ref{eq1})] are one-column matrices. \citet{DenTan99} looked at the absolute sums of the columns of $\mathbf{M}$,
which were termed $J$-characteristics by \citet{TanDen99}.
Specifically, for a resolution $R$ design, these authors introduced $\mathrm{GR}$
as
%
\begin{equation}
\mathrm{GR} = R + 1 - \frac{\max J_{R}}{N},\label{eq3}
\end{equation}
where $J_{R} = |\mathbf{1}_{N}^{\mathrm{T}} \mathbf{M}_{R}|$ is the row
vector of the $J\mbox{-}\mathrm{characteristics}\ |\mathbf{1}_{N}^{\mathrm{T}}
\mathbf{X}_{u_1,\ldots,u_R}|$ obtained from the $ { n
\choose R}$ $R$-factor interaction model columns
$\mathbf{X}_{u_1,\ldots, u_R}$. For 2-level designs, it is straightforward to
verify the following identities:
%
\begin{eqnarray}\label{eq4}
\mathrm{GR} &=& R + 1 - \sqrt{\max_{(u_{1},\ldots,u_{R})}a_{R}(u_{1},
\ldots,u_{R})}
\nonumber
\\[-8pt]
\\[-8pt]
\nonumber
& = & R + 1 - \max_{(u_{1},\ldots,u_{R})}\bigl\llvert \rho
(X_{u_{1}},X_{u_{2},\ldots, u_{R}}) \bigr\rrvert,
\end{eqnarray}
where $\rho$ denotes the correlation; note that the correlation in (\ref{eq4}) does
not depend on which of the $u_{\mathrm{i}}$ takes the role of $u_{1}$. Deng
and Tang [(\citeyear{DenTan99}), Proposition~2] proved a very convincing projection interpretation
of their $\mathrm{GR}$. Unfortunately, Proposition~4.4 of \citet{DieBed02}, in which a particular $\operatorname{OA}(18, 3^{3}, 2)$ is proven to be
indecomposable into two $\operatorname{OA}(9, 3^{3}, 2)$, implies that Deng and Tang's
result cannot be generalized to more than two levels.

The quantitative approach by Evangelaras et al. [(\citeyear{Evaetal05}), their equation (4)]
generalized the correlation version of (\ref{eq4}) by applying it to single df
contrasts for the quantitative factors. For the qualitative factors
considered here, any approach based on direct usage of single df contrasts
is not acceptable because it is coding dependent. The approach for
qualitative factors taken by Evangelaras et al. is unreasonable, as will be
demonstrated in Example~\ref{ex5}. \citet{PanLiu10} also proposed a generalized
resolution based on complex contrasts. For designs with more than 3 levels,
permuting levels for one or more factors will lead to different generalized
resolutions according to their definition, which is unacceptable for
qualitative factors. For 2-level designs, their approach boils down to
omitting the square root from
$\sqrt{\max_{(u_{1},\ldots,u_{R})}a_{R}(u_{1},\ldots,u_{R})}$ in (\ref{eq4}), which
implies that their proposal does not simplify to the well-grounded
generalized resolution of \citet{DenTan99}/\citet{TanDen99} for
2-level designs. This in itself makes their approach unconvincing. Example~\ref{ex5} will compare their approach to ours for 3-level designs. The results from
the previous section can be used to create two adequate generalizations of
$\mathrm{GR}$ for qualitative factors. These are introduced in the following two
definitions.

For the first definition, an $R$ factor projection is considered as
completely aliased, whenever all the levels of at least one of the factors
are completely determined by the level combination of the other $R -1$
factors. Thus, generalized resolution should be equal to $R$, if and only
if there is at least one $R$ factor projection with
$a_{R}(u_{1},\ldots,u_{R}) = s_{\mathrm{min}}-1$. The $\mathrm{GR}$ defined in
Definition~\ref{de2} guarantees this behavior and fulfills all requirements stated
in the \hyperref[sec1]{Introduction}:

\begin{definition}\label{de2}
For an $\operatorname{OA}(N, s_{1}, \ldots, s_{n}, R -1)$,
\[
\mathrm{GR} = R + 1 - \sqrt{\mathop{\max_{\{ u_{1},\ldots,u_{R}\} }}_{
 \subseteq \{ 1,\ldots,n\}} \frac{a_{R} ( u_{1},\ldots,u_{R}
)}{\mathop{\mathrm{min}}\limits_{i = 1,\ldots,R}s_{u_{i}} - 1}}.
\]

In words, $\mathrm{GR}$ increases the resolution by one minus the square root of the
worst case average $R^{2}$ obtained from any $R$ factor projection, when
regressing the main effects columns in orthogonal coding from a factor with
the minimum number of levels on the other factors in the projection. It is
straightforward to see that (\ref{eq4}) is a special case of the definition, since
the denominator is 1 for 2-level designs. Regarding the requirements stated
in the \hyperref[sec1]{Introduction}, (i) $\mathrm{GR}$ from Definition~\ref{de2} is coding invariant because
the $a_{R}(\cdot)$ are coding invariant according to \citet{XuWu01}. (ii) The
technique is obviously applicable for symmetric and asymmetric designs
alike, and (iii) $\mathrm{GR} < R + 1$ follows from the resolution, $\mathrm{GR} \geq R$
follows from part (ii) of Corollary~\ref{co1}, $\mathrm{GR} = R$ is equivalent to complete
confounding in at least one $R$-factor projection according to part (i) of
Corollary~\ref{co1}.
\end{definition}

\setcounter{example}{3}
\begin{example}\label{ex4}
The $\mathrm{GR}$ values for the designs from Examples \ref{ex1} and
\ref{ex3} are $3 (\mathrm{GR}=R)$, the $\mathrm{GR}$ value for the design from Example~\ref{ex2} is $3 + 1 -
\sqrt{4 / 9} = 3.33$, and the $\mathrm{GR}$ values for all designs from Table~\ref{tab3} are
$3 + 1 - \sqrt{1 / 3} = 3.42$.
\end{example}

Now, complete aliasing is considered regarding individual degrees of
freedom (df). A coding invariant individual df approach considers a
factor's main effect as completely aliased in an $R$ factor projection,
whenever there is at least one pair of canonical variates with correlation
one. A projection is considered completely aliased, if at least one
factor's main effect is completely aliased in this individual df sense.
Note that it is now possible that factors with the same number of levels
can show different extents of individual df aliasing within the same
projection, as will be seen in Example~\ref{ex5} below.

\begin{definition}\label{de3}
For an $\operatorname{OA}(N, s_{1}, \ldots, s_{n}, R -1)$ and
tuples $(c, \mathrm{C})$ with $\mathrm{C} = \{u_{1},\ldots,u_{R}\} \setminus
\{c\}$,
\[
\mathrm{GR}_{\mathrm{ind}} = R + 1 - \max_{\{ u_{1},\ldots,u_{R}\} \subseteq
\{ 1,\ldots,n\}} \max
_{c \in \{ u_{1},\ldots,u_{R}\}} r_{1} ( \mathbf{X}_{c};
\mathbf{F}_{\mathrm{C}} )
\]
with $r_{1}(\mathbf{X}_{c}; \mathbf{F}_{\mathrm{C}})$ the largest canonical
correlation between the main effects model matrix for factor $c$ and the
full model matrix of the factors in $\mathrm{C}$.
\end{definition}

In words, $\mathrm{GR}_{\mathrm{ind}}$ is the worst case confounding for an individual
main effects df in the design that can be obtained by the worst case coding
(which corresponds to the $\mathbf{v}_{1}$ vector associated with the worst
canonical correlation). Obviously, $\mathrm{GR}_{\mathrm{ind}}$ is thus a stricter
criterion than $\mathrm{GR}$. Formally, Theorem~\ref{th3} implies that $\mathrm{GR}$ from Definition~\ref{de2} can be written as
%
\begin{equation}
\mathrm{GR} = R + 1 - \sqrt{\mathop{\max_{( u_{1},\ldots,u_{R}  )\dvtx
}}_{\{ u_{1},\ldots,u_{R}\} \subseteq \{ 1,\ldots,n\}} \frac{\sum_{j = 1}^{s_{u_{1}} - 1} r_{j} (
\mathbf{X}_{u_{1}};\mathbf{F}_{ \{ u_{2},\ldots,u_{R}  \}}
)^{2}} {\mathrm{min}_{i}s_{u_{i}} - 1}}.
\label{eq5}
\end{equation}

Note that maximization in (\ref{eq5}) is over tuples, so that it is ensured that
the factor with the minimum number of levels does also get into the first
position. Comparing~(\ref{eq5}) with Definition~\ref{de3}, $\mathrm{GR}_{\mathrm{ind}}\leq \mathrm{GR}$ is
obvious, because $r_{1}^{2}$ cannot be smaller than the average over all
$r_{i}^{2}$ (but can be equal, if all canonical correlations have the same
size). This is stated in a theorem.

\begin{theorem}\label{th4}
For $\mathrm{GR}$ from Definition~\ref{de2} and $\mathrm{GR}_{\mathrm{ind}}$ from
Definition~\ref{de3}, $\mathrm{GR}_{\mathrm{ind}} \leq \mathrm{GR}$.
\end{theorem}

\begin{remark}\label{re3}
(i) Under normalized orthogonal coding, the full
model matrix $\mathbf{F}_{\mathrm{C}}$ in Definition~\ref{de3} can again be
replaced by the $R -1$ factor interaction matrix~$\mathbf{X}_{\mathrm{C}}$.
(ii) Definition~\ref{de3} involves calculation of
$R  {n \choose
R}$ canonical correlations ($R$ correlations for each $R$
factor projection). In any projection with at least one 2-level factor, it
is sufficient to calculate one single canonical correlation obtained with
an arbitrary 2-level factor in the role of $Y$, because this is necessarily
the worst case. Nevertheless, calculation of $\mathrm{GR}_{\mathrm{ind}}$ carries some
computational burden for designs with many factors.

\begin{table}[b]
\tabcolsep=0pt
\caption{The Taguchi L18 (transposed)}
\label{tab4}
\begin{tabular*}{\textwidth}{@{\extracolsep{\fill}}lcccccccccccccccccc@{}}
\hline
&\multicolumn{18}{c@{}}{\textbf{Row}}\\[-6pt]
&\multicolumn{18}{c@{}}{\hrulefill}\\
\textbf{Column}&
\textbf{1} & \textbf{2} & \textbf{3} & \textbf{4} & \textbf{5} & \textbf{6} &
\textbf{7} & \textbf{8} & \textbf{9} & \textbf{10} & \textbf{11} &
\textbf{12} & \textbf{13} & \textbf{14} & \textbf{15} & \textbf{16} &
\textbf{17} & \textbf{18}\\
\hline
1 & 0 & 0 & 0 & 0 & 0 & 0 & 0 & 0 & 0 & 1 & 1 & 1 & 1 & 1 & 1 & 1 & 1 & 1\\
2 & 0 & 0 & 0 & 1 & 1 & 1 & 2 & 2 & 2 & 0 & 0 & 0 & 1 & 1 & 1 & 2 & 2 & 2\\
3 & 0 & 1 & 2 & 0 & 1 & 2 & 0 & 1 & 2 & 0 & 1 & 2 & 0 & 1 & 2 & 0 & 1 & 2\\
4 & 0 & 1 & 2 & 0 & 1 & 2 & 1 & 2 & 0 & 2 & 0 & 1 & 1 & 2 & 0 & 2 & 0 & 1\\
5 & 0 & 1 & 2 & 1 & 2 & 0 & 0 & 1 & 2 & 2 & 0 & 1 & 2 & 0 & 1 & 1 & 2 & 0\\
6 & 0 & 1 & 2 & 1 & 2 & 0 & 2 & 0 & 1 & 1 & 2 & 0 & 0 & 1 & 2 & 2 & 0 & 1\\
7 & 0 & 1 & 2 & 2 & 0 & 1 & 1 & 2 & 0 & 1 & 2 & 0 & 2 & 0 & 1 & 0 & 1 & 2\\
8 & 0 & 1 & 2 & 2 & 0 & 1 & 2 & 0 & 1 & 0 & 1 & 2 & 1 & 2 & 0 & 1 & 2 & 0\\
\hline
\end{tabular*}
\end{table}

Obviously, (\ref{eq4}) is a special case of $\mathrm{GR}_{\mathrm{ind}}$, since the average
$R^{2}$ coincides with the only squared canonical correlation for
projections of $R$ $2$-level factors. $\mathrm{GR}_{\mathrm{ind}}$ also fulfills all
requirements stated in the \hyperref[sec1]{Introduction}: (i) $\mathrm{GR}_{\mathrm{ind}}$ is coding
invariant because the canonical correlations are invariant to affine
transformations of the $X$ and $Y$ variables, as was discussed in Section~\ref{sec2}. (ii) The technique is obviously applicable for symmetric and asymmetric
designs alike, and (iii) $\mathrm{GR}_{\mathrm{ind}} < R+1$ again follows from the
resolution, $\mathrm{GR}_{\mathrm{ind}} \geq R$ follows from the properties of
correlations, and $\mathrm{GR}_{\mathrm{ind}} = R$ is obviously equivalent to complete
confounding of at least one main effects contrast in at least one $R$
factor projection, in the individual df sense discussed above.
\end{remark}

\begin{example}\label{ex5}
We consider the three nonisomorphic $\operatorname{OA}(18, 3^{3},
2)$ that can be obtained as projections from the well-known Taguchi L18 (see
Table~\ref{tab4}) by using columns 3, 4 and 5 $(D_{1})$, columns 2, 3 and 6
$(D_{2})$ or columns 2, 4 and 5\vadjust{\goodbreak} $(D_{3})$. We have $A_{3}(D_{1}) = 0.5$,
$A_{3}(D_{2}) = 1$ and $A_{3}(D_{3}) = 2$, and consequently $\mathrm{GR}(D_{1}) =
3.5$, $\mathrm{GR}(D_{2}) = 3.29$ and $\mathrm{GR}(D_{3}) = 3$. For calculating
$\mathrm{GR}_{\mathrm{ind}}$, the largest canonical correlations of all factors in the
role of $Y$ are needed. These are all 0.5 for $D_{1}$ and all 1 for
$D_{3}$, such that $\mathrm{GR}_{\mathrm{ind}} = \mathrm{GR}$ for these two designs. For
$D_{2}$, the largest canonical correlation is 1 with the first factor (from
column 2 of the L18) in the role of $Y$, while it is $\sqrt{0.5}$ with
either of the other two factors in the role of $Y$; thus, $\mathrm{GR}_{\mathrm{ind}} =
3 < \mathrm{GR} = 3.29$. The completely aliased 1 df contrast of the first factor is
the contrast of the third level vs. the other two levels, which is apparent
from Table~\ref{tab5}: the contrast $\mathrm{A} = 2$ vs. $\mathrm{A}$ in $(0,1)$ is fully aliased
with the contrast of one level of B vs. the other two, given a particular
level of C. Regardless of factor coding, this direct aliasing is reflected
by a canonical correlation ``one'' for the first canonical variate of the
main effects contrast matrix of factor A.
\end{example}

\begin{table}
\caption{Frequency table of columns $2\ (=\mathrm{A})$, $3\ (=\mathrm{B})$ and $6\ (=\mathrm{C})$ of the
Taguchi L18}\label{tab5}
\begin{tabular*}{\textwidth}{@{\extracolsep{\fill}}ccc@{}}
\textbf{, ,} $\mathbf{C \bolds{=} 0}$&\textbf{, ,} $\mathbf{C \bolds{=} 1}$&
\textbf{, ,} $\mathbf{C \bolds{=} 2}$\\
\begin{tabular}{ccccc}
&\textbf{B}&&&\\
\textbf{A}&& \textbf{0} &\textbf{1} &\textbf{2}\\
&\textbf{0}& 1& 0 &1\\
&\textbf{1}& 1 &0 &1\\
&\textbf{2}& 0& 2& 0
\end{tabular}
&
\begin{tabular}{ccccc}
&\textbf{B}&&&\\
\textbf{A}&& \textbf{0} &\textbf{1} &\textbf{2}\\
&\textbf{0}& 1& 1 &0\\
&\textbf{1}& 1 &1&0\\
&\textbf{2}& 0& 0& 2
\end{tabular}
&
\begin{tabular}{ccccc}
&\textbf{B}&&&\\
\textbf{A}&& \textbf{0} &\textbf{1} &\textbf{2}\\
&\textbf{0}& 0& 1 &1\\
&\textbf{1}& 0 &1 &1\\
&\textbf{2}& 2& 0& 0
\end{tabular}
\end{tabular*}
\end{table}

Using this example, we now compare the $\mathrm{GR}$ introduced here to proposals by
\citet{Evaetal05} and \citet{PanLiu10}: The \textit{GRes} values
reported by \citet{Evaetal05} for designs $D_{1}$, $D_{2}$ and
$D_{3}$ in the qualitative case are 3.75, 3.6464, 3.5, respectively;
especially the 3.5 for the completely aliased design $D_{3}$ does not make
sense. Pang and Liu reported values 3.75, 3.75 and 3, respectively; here,
at least the completely aliased design $D_{3}$ is assigned the value~``3.''
Introducing the square root, as was discussed in connection with equation
(\ref{eq4}), their generalized resolutions become 3.5, 3.5 and 3, respectively,
that is, they coincide with our $\mathrm{GR}$ results for designs $D_{1}$ and
$D_{3}$. For design $D_{2}$, their value 3.5 is still different from our
3.29 for the following reason: our approach considers $A_{3} =
a_{3}(1,2,3)$ as a sum of two $R^{2}$-values and subtracts the square root
of their average or maximum ($\mathrm{GR}$ or $\mathrm{GR}_{\mathrm{ind}}$, resp.), while
Pang and Liu's approach considers it as a sum of $2^{3} =8$ summands,
reflecting the potentially different linear combinations of the three
factors in the Galois field sense, the (square root of the) maximum of
which they subtract from $R+1$.

\section{Properties of GR}\label{sec4}

Let G be the set of all runs of an $s_{1} \times \cdots \times s_{n}$ full
factorial design, with $\llvert  \mathrm{G} \rrvert  = \prod_{i = 1}^{n}
s_{i}$ the cardinality of G. For any design $D$ in $N$ runs for $n$ factors
at $s_{1}, \ldots, s_{n}$ levels, let $N_{\mathbf{x}}$ be the number of
times that a point $\mathbf{x} \in\mathrm{G}$ appears in $D$. $\bar{N} = N /
\llvert  \mathrm{G} \rrvert $ denotes the average frequency for each point of
G in the design $D$. We can measure the goodness of a fractional factorial
design $D$ by the uniformity of the design points of $D$ in the set of all
points in G, that is, the uniformity of the frequency distribution
$N_{\mathbf{x}}$. One measure, suggested by \citet{Tan01} and \citet{AiZha04}, is the variance
\[
\operatorname{V} (D) = \frac{1}{\llvert  \mathrm{G} \rrvert }\sum_{\mathbf{x}
\in \mathrm{G}}
( N_{\mathbf{x}} - \bar{N} )^{2} = \frac{1}{\llvert  \mathrm{G} \rrvert }\sum
_{\mathbf{x} \in \mathrm{G}} N_{\mathbf{x}}^{2} - \bar{N}^{2}.
\]

Let $N = q |\mathrm{G}| + r$ with nonnegative integer $q$ and $r$ and $0
\leq r < |G|$ (often $q=0$), that is, $r = N \operatorname{mod} |\mathrm{G}|$ is the remainder of
$N$ when divided by $|\mathrm{G}|$. Note that $\sum_{\mathbf{x} \in \mathrm{G}}
N_{\mathbf{x}} = N$, so $\mathrm{V}(D)$ is minimized if and only if each
$N_{\mathbf{x}}$ takes values on $q$ or $q + 1$ for any $\mathbf{x} \in \mathrm{G}$.
When $r$ points in $\mathrm{G}$ appear $q + 1$ times and the remaining $|\mathrm{G}| -
r$ points appear $q$ times, $\mathrm{V}(D)$ reaches the minimal value
$r(|\mathrm{G}| - r) / |\mathrm{G}|^{2}$. \citet{AiZha04} showed that
$\mathrm{V}(D)$ is a function of GWLP. In particular, if $D$ has strength
$n - 1$, their result implies that $V(D) = \bar{N}^{2}A_{n}(d)$. Combining
these results, and using the following definition, we obtain an upper bound
for $\mathrm{GR}$ for some classes of designs and provide a necessary and sufficient
condition under which this bound is achieved.

\begin{definition}[{[Modified from \citet{Xu03}]}]\label{de4}
(i) A design $D$ has maximum $t$-balance, if and only if the
possible level combinations for all projections onto $t$ columns occur as
equally often as possible, that is, either $q$ or $q+1$ times, where $q$ is
an integer such that $N = q|\mathrm{G}_{\mathrm{proj}}| + r$ with $\mathrm{G}_{\mathrm{proj}}$ the set
of all runs for the full factorial design of each respective
$t$-factor-projection and $0 \leq r < | \mathrm{G}_{\mathrm{proj}}|$.

(ii) An $\operatorname{OA}(N, s_{1}, \ldots, s_{n}, t -1)$ with $n \geq t$ has
weak strength $t$ if and only if it has maximum $t$-balance. We denote weak
strength $t$ as $\operatorname{OA}(N,s_{1}, \ldots, s_{n}, t^{-}$).
\end{definition}

\begin{remark}\label{re4}
\citet{Xu03} did not require strength
${t}-1$ in the definition of weak strength
${t}$, that is, the \citet{Xu03} definition of
weak strength ${t}$ corresponds to our
definition of maximum ${t}$-balance. For the
frequent case, for which all ${t}$-factor
projections have ${q}=0$ or
${q}=1$ and
${r}=0$ in Definition~\ref{de4}(i), maximum
${t}$-balance is equivalent to the absence of
repeated runs in any projection onto ${t}$
factors. In that case, maximum
${t}$-balance implies maximum
${k}$-balance for ${k}
> {t}$, and weak strength
${t}$ is equivalent to strength
${t}-1$ with absence of repeated runs in any
projection onto ${t}$ or more factors.
\end{remark}

\begin{theorem}\label{th5}
Let $D$ be an $\operatorname{OA}(N, s_{1}, \ldots, s_{R}, R - 1)$.
Then $A_{R} ( D  ) \ge \frac{r ( \prod_{i = 1}^{R} s_{i} - r
)}{N^{2}}$, where $r$ is the remainder when $N$ is divided by
$\prod_{i = 1}^{R} s_{i}$. The equality holds if and only if $D$ has weak
strength $R$.
\end{theorem}

As all $R$ factor projections of any $\operatorname{OA}(N, s_{1}, \ldots, s_{n}, R^{-})$
fulfill the necessary and sufficient condition of Theorem~\ref{th5}, we have the
following corollary.\vadjust{\goodbreak}

\begin{corollary}\label{co2}
Suppose that an $\operatorname{OA}(N, s_{1},\ldots,  s_{n}, R)$
does not exist. Then any $\operatorname{OA}(N, s_{1}, \ldots, s_{n}, R^{-})$ has maximum
$\mathrm{GR}$ among all $\operatorname{OA}(N, s_{1},\ldots, s_{n},\break  R - 1)$.
\end{corollary}

\begin{corollary}\label{co3}
Suppose that an $\operatorname{OA}(N, s^{n}, R)$ does not
exist. Let $D$ be an $\operatorname{OA}(N, s^{n}, R -1)$. Then $\mathrm{GR}(D) \le R + 1 -
\sqrt{\frac{r ( s^{R} - r  )}{N^{2} ( s - 1  )}}$, where
$r$ is the remainder when $N$ is divided by $s^{R}$. The equality holds if
and only if $D$ has weak strength $R$.
\end{corollary}

\begin{example}\label{ex6}
(1) Any projection onto three 3-level columns from
an $\operatorname{OA}(18, 6^{1}3^{6}, 2)$ has 18 distinct runs
($q=0$, $r=N=18$) and is an
{OA} of weak strength 3, so it has $A_{3} =1/2$ and $\mathrm{GR} =
 4 - \sqrt{18
\cdot 9 /  ( 18^{2} \cdot 2  )}=3.5$. (2) Any projection onto
three or more $s$-level columns from an $\operatorname{OA}(s^{2}, s^{s+1}, 2)$ has $\mathrm{GR} =
3$, since $N = r = s^{2}$, so that the upper limit from the corollary
becomes $\mathrm{GR} = R = 3$.
\end{example}

Using the following lemma according to \citet{MukWu95}, Corollary~\ref{co3}
can be applied to a further class of designs.

\begin{lemma}[{[\citet{MukWu95}]}]\label{le1} For a saturated $\operatorname{OA}(N,
s_{1}^{n_1}\mathrm{s}_{2}^{n_2}, 2)$ with $n_{1}(s_{1} -1) + n_{2}(s_{2}
-1) = N - 1$, let $\delta_{i}(a, b)$ be the number of coincidences of two
distinct rows $a$ and $b$ in the $n_{i}$ columns of $s_{i}$ levels, for $i
= 1, 2$. Then
\[
s_{1} \delta_{1}(a, b) + s_{2}
\delta_{2}(a, b) = n_{1} + n_{2} - 1.
\]
\end{lemma}

Consider a saturated $\operatorname{OA}(2s^{2}, (2s)^{1} s^{2s}, 2)$, where $r = N =
2s^{2}, s_{1} = 2s, s_{2} = s, n_{1} = 1, n_{2} = 2s$. From Lemma~\ref{le1}, we have $2\delta_{1}(a, b) + \delta_{2}(a, b) = 2$. So any
projection onto three or more $s$-level columns has no repeated runs, and
thus it achieves the upper limit $\mathrm{GR} = 4 - \sqrt{ ( s - 2  ) /
( 2s - 2  )}$ according to Corollary~\ref{co3}.

\begin{corollary}\label{co4}
For a saturated $\operatorname{OA}(2s^{2}, (2s)^{1} s^{2s},
2)$, any projection onto three or more $s$-level columns has $\mathrm{GR} = 4 -
\sqrt{ ( s - 2  ) /  ( 2s - 2  )}$, which is optimum
among all possible OAs in $2s^{2}$ runs.
\end{corollary}

\begin{example}\label{ex7}
Design 1 of Table~\ref{tab3} is isomorphic to a projection
from a saturated $\operatorname{OA}(32,
8^{1}4^{8}, 2)$.
${A}_{3}$ attains the lower bound
from Theorem \ref{th5} $(32\cdot (64-32)/32^{2}= 1)$,
and thus GR attains the upper bound
$4-(1/3)^{\textonehalf} = 3.42$ from the
corollary.
\end{example}

Because of Theorem~\ref{th4}, any upper bound for $\mathrm{GR}$ is of course also an upper
bound for $\mathrm{GR}_{\mathrm{ind}}$, that is, Corollaries \ref{co3} and \ref{co4} also provide upper
bounds for $\mathrm{GR}_{\mathrm{ind}}$. However, for $\mathrm{GR}_{\mathrm{ind}}$ the bounds are not
tight in general; for example, $\mathrm{GR}_{\mathrm{ind}} = 3$ for the design of
Example~\ref{ex7} (see also Example~\ref{ex9} in the following section).

\citet{But05} previously showed that all projections onto $s$-level columns
of $\operatorname{OA}(s^{2}, s^{s+1}, 2)$ or $\operatorname{OA}(2s^{2}, (2s)^{1} s^{2s}, 2)$ have GMA
among all possible designs.

\section{Factor wise GR values}\label{sec5}

In Section~\ref{sec3}, two versions of overall generalized resolution were defined:
$\mathrm{GR}$ and $\mathrm{GR}_{\mathrm{ind}}$. These take a worst case perspective: even if a
single projection in a large design is completely confounded---in the case
of mixed level designs or $\mathrm{GR}_{\mathrm{ind}}$ affecting perhaps only one factor
within that projection---the overall metric takes the worst case value
$R$. It can therefore be useful to accompany $\mathrm{GR}$ and $\mathrm{GR}_{\mathrm{ind}}$ by
factor specific summaries. For the factor specific individual df
perspective, one simply has to omit the maximization over the factors in
each projection and has to use the factor of interest in the role of $Y$
only. For a factor specific complete confounding perspective, one has to
divide each projection's $a_{R}(\cdot)$ value by the factor's df rather than the
minimum df, in order to obtain the average $R^{2}$ value for this
particular factor. This leads to the following definition.

\begin{definition}\label{de5}
For an $\operatorname{OA}(N, s_{1},\ldots, s_{n}, R -1)$,
define
\begin{longlist}[(ii)]
\item[(i)] $\mathrm{GR}_{\mathrm{tot}(i)} = R + 1 - \sqrt{\max_{  \{
u_{2},\ldots,u_{R}  \}  \subseteq  \{ 1,\ldots,n
\} \setminus  \{ \mathrm{i}  \}}
\frac{a_{R} (
i,u_{2},\ldots,u_{R}  )}{s_{i} - 1}}$,

\item[(ii)] $\mathrm{GR}_{\mathrm{ind}(i)} = R + 1 - \max_{\{ i,u_{2},\ldots,u_{R}\}
\subseteq \{ 1,\ldots,n\}} r_{1} (
\mathbf{X}_{i};\mathbf{X}_{u_{2},\ldots, u_{R}}  )$, with $\mathbf{X}_{i}$
the model matrix of factor $i$ and $\mathbf{X}_{u_2,\ldots, u_R}$ the $R -1$
factor interaction model matrix of the factors in $\{u_{2},\ldots,u_{R}\}$
in normalized orthogonal coding, and $r_{1}(\mathbf{Y};\mathbf{X})$ the
first canonical correlation between matrices $\mathbf{X}$ and
$\mathbf{Y}$.
\end{longlist}
\end{definition}

It is straightforward to verify that $\mathrm{GR}$ and $\mathrm{GR}_{\mathrm{ind}}$ can be
calculated as the respective minima of the factor specific $\mathrm{GR}$ values from
Definition~\ref{de5}.

\begin{theorem}\label{th6}
For the quantities from Definitions \ref{de2}, \ref{de3} and \ref{de5}, we
have
\begin{longlist}[(ii)]
\item[(i)] $\mathrm{GR} = \mathrm{min}_{i} \mathrm{GR}_{\mathrm{tot}(i)}$,
\item[(ii)] $\mathrm{GR}_{\mathrm{ind}} = \mathrm{min}_{i} \mathrm{GR}_{\mathrm{ind}(i)}$.
\end{longlist}
\end{theorem}

\begin{example}\label{ex8}
The Taguchi L18 has $\mathrm{GR} = \mathrm{GR}_{\mathrm{ind}} = 3$, and
the following $\mathrm{GR}_{\mathrm{ind}(i)}$ and $\mathrm{GR}_{\mathrm{tot}(i)}$ values
($\mathrm{GR}_{\mathrm{ind}(i)} = \mathrm{GR}_{\mathrm{tot}(i)}$ for all $i$): 3.18, 3, 3.29, 3,
3, 3.29, 3.29, 3.29. When omitting the second column, the remaining seven
columns have $\mathrm{GR} = \mathrm{GR}_{\mathrm{ind}} = 3.18$, again with $\mathrm{GR}_{\mathrm{ind}(i)}
= \mathrm{GR}_{\mathrm{tot}(i)}$ and the value for all 3-level factors at 3.42. When
omitting the fourth column instead, the then remaining seven columns have
$\mathrm{GR} = 3.18$, $\mathrm{GR}_{\mathrm{ind}} = 3$, $\mathrm{GR}_{\mathrm{tot}(i)}$ values 3.18, 3.29,
3.29, 3.42, 3.29, 3.29, 3.29 and $\mathrm{GR}_{\mathrm{ind}(i)}$ values the same,
except for the second column, which has $\mathrm{GR}_{\mathrm{ind}(2)} = 3$.
\end{example}

\begin{table}
\caption{Largest canonical correlations, $\mathrm{GR}_{\mathrm{ind}(i)}$ and
$\mathrm{GR}_{\mathrm{ind}}$ values for the GMA
$\operatorname{OA}(32, 4^{3}, 2)$}
\label{tab6}
\begin{tabular*}{\textwidth}{@{\extracolsep{\fill}}lccccccc@{}}
\hline
& $\bolds{r_{1}(1;23)}$ & $\bolds{r_{1}(2;13)}$ &
$\bolds{r_{1}(3;12)}$ & $\bolds{\mathrm{GR}_{\mathrm{ind}(1)}}$ &
$\bolds{\mathrm{GR}_{\mathrm{ind}(2)}}$ & $\bolds{\mathrm{GR}_{\mathrm{ind}(3)}}$ &
 $\bolds{\mathrm{GR}_{\mathrm{ind}}}$\\
 \hline
\phantom{0}1 & 1.000 & 1.000 & 1.000 & 3.000 & 3.000 & 3.000 & 3.000\\
\phantom{0}2 & 0.866 & 0.866 & 0.866 & 3.134 & 3.134 & 3.134 & 3.134\\
\phantom{0}3 & 0.707 & 0.707 & 1.000 & 3.293 & 3.293 & 3.000 & 3.000\\
\phantom{0}4 & 0.707 & 0.707 & 0.866 & 3.293 & 3.293 & 3.134 & 3.134\\
\phantom{0}5 & 0.707 & 0.707 & 0.791 & 3.293 & 3.293 & 3.209 & 3.209\\
\phantom{0}6 & 0.707 & 0.707 & 0.707 & 3.293 & 3.293 & 3.293 & 3.293\\
\phantom{0}7 & 0.707 & 0.707 & 0.707 & 3.293 & 3.293 & 3.293 & 3.293\\
\phantom{0}8 & 0.707 & 0.707 & 0.707 & 3.293 & 3.293 & 3.293 & 3.293\\
\phantom{0}9 & 0.612 & 0.612 & 0.612 & 3.388 & 3.388 & 3.388 & 3.388\\
10 & 0.707 & 0.707 & 0.707 & 3.293 & 3.293 & 3.293 & 3.293\\
\hline
\end{tabular*}
\end{table}

$\mathrm{GR}$ from Definition~\ref{de2} and $\mathrm{GR}_{\mathrm{ind}}$ from Definition~\ref{de3} are not the
only possible generalizations of (\ref{eq4}). It is also possible to define a
$\mathrm{GR}_{\mathrm{tot}}$, by declaring only those $R$ factor projections as
completely confounded for which all factors are completely confounded. For
this, the factor wise average $R^{2}$ values for each projection---also
used in $\mathrm{GR}_{\mathrm{tot}(i)}$---need to be considered. A projection is
completely confounded, if these are all one, which can be formalized by
requesting their minimum or their average to be one. The average appears
more informative, leading to
%
\begin{equation}
\mathrm{GR}_{\mathrm{tot}} = R + 1 - \sqrt{\mathop{\max_{\{
u_{1},\ldots,u_{R} \}}}_{  \subseteq \{ 1,\ldots,n
\}}
\frac{1}{R}\sum_{i = 1}^{R}
\frac{a_{R}( u_{1},\ldots,u_{R}
)}{s_{{u}_{i}} - 1}}. \label{eq6}
\end{equation}

It is straightforward to see that $\mathrm{GR}_{\mathrm{tot}} \geq \mathrm{GR}$, and that
$\mathrm{GR}_{\mathrm{tot}} = \mathrm{GR}$ for symmetric designs. The asymmetric design of Table~\ref{tab2}
(Example~\ref{ex3}) has $\mathrm{GR} = 3$ and $\mathrm{GR}_{\mathrm{tot}} = 3 + 1 - \sqrt{ ( 1 + 1 + 1 /
3  ) / 3} = 3.12> 3$, in spite of the fact that two of its factors are
completely confounded. Of course, mixed level projections can never be
completely confounded according to (\ref{eq6}), which is the main reason why we
have not pursued this approach.

The final example uses the designs of Table~\ref{tab3} to show that $\mathrm{GR}_{\mathrm{ind}}$
and the $\mathrm{GR}_{\mathrm{ind}(i)}$ can introduce meaningful differentiation
between GMA designs.

\begin{example}\label{ex9}
All designs of Table~\ref{tab3} had $A_{3} =1$ and $\mathrm{GR}=3.42$.
The information provided in Table~\ref{tab3} is insufficient for determining
$\mathrm{GR}_{\mathrm{ind}}$. Table~\ref{tab6} provides the necessary information: the largest
canonical correlations are the same regardless which variable is chosen as
the $Y$ variable for seven designs, while they vary with the choice of the
$Y$ variable for three designs. There are five different $\mathrm{GR}_{\mathrm{ind}}$
values for these 10 designs that were not further differentiated by $A_{3}$
or $\mathrm{GR}$, and in combination with the $\mathrm{GR}_{\mathrm{ind}(i)}$, seven different
structures can be distinguished.
\end{example}

The differentiation achieved by $\mathrm{GR}_{\mathrm{ind}}$ is meaningful, as can be
seen by comparing frequency tables of the first, third and ninth design
(see Table~\ref{tab7}). The first and third design have $\mathrm{GR}_{\mathrm{ind}}=3$, which is
due to a very regular confounding pattern: in the first design,
dichotomizing each factor into a $0/1$ vs. $2/3$ design yields a regular
resolution III 2-level design (four different runs only), that is, each
main effect contrast $0/1$ vs. $2/3$ is completely confounded by the two-factor
interaction of the other two $0/1$ vs. $2/3$ contrasts; the third design shows
this severe confounding for factor C only, whose $0/1$ vs. $2/3$ contrast
is likewise completely confounded by the interaction between factors A and
B. Design 9 is the best of all GMA designs in terms of $\mathrm{GR}_{\mathrm{ind}}$. It
does not display such a strong regularity in behavior. $\mathrm{GR}_{\mathrm{ind}}$
treats designs 1 and 3 alike, although design 1 is clearly more severely
affected than design 3, which can be seen from the individual
$\mathrm{GR}_{\mathrm{ind}(i)}$. However, as generalized resolution has always taken a
``worst case'' perspective, this way of handling things is appropriate in
this context.

\begin{table}\tabcolsep=4pt
\caption{Frequency tables of designs 1, 3 and 9
from Table~\protect\ref{tab6}}\label{tab7}
\begin{tabular*}{\textwidth}{@{\extracolsep{\fill}}cccc@{}}
\textbf{, ,} $\mathbf{C \bolds{=} 0}$&\textbf{, ,} $\mathbf{C \bolds{=} 1}$&
\textbf{, ,} $\mathbf{C \bolds{=} 2}$&\textbf{, ,} $\mathbf{C \bolds{=} 3}$\\
\multicolumn{4}{@{}l}{Design 1}\\
\begin{tabular}{@{}cccccc@{}}
&\textbf{B}&&&&\\
\textbf{A}&& \textbf{0} &\textbf{1} &\textbf{2}&\textbf{3}\\
&\textbf{0}& {1}& {1} &{0}&{0}\\
&\textbf{1}& {1} &{1} &{0}&{0}\\
&\textbf{2}& {0}& {0}& {1}&{1}\\
&\textbf{3}&{0}&{0}&{1}&{1}
\end{tabular}
&
\begin{tabular}{@{}cccccc@{}}
&\textbf{B}&&&&\\
\textbf{A}&& \textbf{0} &\textbf{1} &\textbf{2}&\textbf{3}\\
&\textbf{0}& {1}& {1} &{0}&{0}\\
&\textbf{1}& {1} &{1} &{0}&{0}\\
&\textbf{2}& {0}& {0}& {1}&{1}\\
&\textbf{3}&{0}&{0}&{1}&{1}
\end{tabular}
&
\begin{tabular}{@{}cccccc@{}}
&\textbf{B}&&&&\\
\textbf{A}&& \textbf{0} &\textbf{1} &\textbf{2}&\textbf{3}\\
&\textbf{0}& {0}& {0} &{1}&{1}\\
&\textbf{1}& {0} &{0} &{1}&{1}\\
&\textbf{2}& {1}& {1}& {0}&{0}\\
&\textbf{3}&{1}&{1}&{0}&{0}
\end{tabular}
&
\begin{tabular}{@{}cccccc@{}}
&\textbf{B}&&&&\\
\textbf{A}&& \textbf{0} &\textbf{1} &\textbf{2}&\textbf{3}\\
&\textbf{0}& {0}& {0} &{1}&{1}\\
&\textbf{1}& {0} &{0} &{1}&{1}\\
&\textbf{2}& {1}& {1}& {0}&{0}\\
&\textbf{3}&{1}&{1}&{0}&{0}
\end{tabular}\\[6pt]
\multicolumn{4}{@{}l}{Design 3}\\
\begin{tabular}{@{}cccccc@{}}
&\textbf{B}&&&&\\
\textbf{A}&& \textbf{0} &\textbf{1} &\textbf{2}&\textbf{3}\\
&\textbf{0}& {1}& {1} &{0}&{0}\\
&\textbf{1}& {1} &{0} &{1}&{0}\\
&\textbf{2}& {0}& {1}& {0}&{1}\\
&\textbf{3}&{0}&{0}&{1}&{1}
\end{tabular}
&
\begin{tabular}{@{}cccccc@{}}
&\textbf{B}&&&&\\
\textbf{A}&& \textbf{0} &\textbf{1} &\textbf{2}&\textbf{3}\\
&\textbf{0}& {1}& {1} &{0}&{0}\\
&\textbf{1}& {1} &{0} &{1}&{0}\\
&\textbf{2}& {0}& {1}& {0}&{1}\\
&\textbf{3}&{0}&{0}&{1}&{1}
\end{tabular}&
\begin{tabular}{@{}cccccc@{}}
&\textbf{B}&&&&\\
\textbf{A}&& \textbf{0} &\textbf{1} &\textbf{2}&\textbf{3}\\
&\textbf{0}& {0}& {0} &{1}&{1}\\
&\textbf{1}& {0} &{1} &{0}&{1}\\
&\textbf{2}& {1}& {0}& {1}&{0}\\
&\textbf{3}&{1}&{1}&{0}&{0}
\end{tabular}&
\begin{tabular}{@{}cccccc@{}}
&\textbf{B}&&&&\\
\textbf{A}&& \textbf{0} &\textbf{1} &\textbf{2}&\textbf{3}\\
&\textbf{0}& {0}& {0} &{1}&{1}\\
&\textbf{1}& {0} &{1} &{0}&{1}\\
&\textbf{2}& {1}& {0}& {1}&{0}\\
&\textbf{3}&{1}&{1}&{0}&{0}
\end{tabular}\\[6pt]
\multicolumn{4}{@{}l}{Design 9}\\
\begin{tabular}{@{}cccccc@{}}
&\textbf{B}&&&&\\
\textbf{A}&& \textbf{0} &\textbf{1} &\textbf{2}&\textbf{3}\\
&\textbf{0}& {1}& {1} &{0}&{0}\\
&\textbf{1}& {1} &{0} &{1}&{0}\\
&\textbf{2}& {0}& {0}& {1}&{1}\\
&\textbf{3}&{0}&{1}&{0}&{1}
\end{tabular}&
\begin{tabular}{@{}cccccc@{}}
&\textbf{B}&&&&\\
\textbf{A}&& \textbf{0} &\textbf{1} &\textbf{2}&\textbf{3}\\
&\textbf{0}& {1}& {0} &{1}&{0}\\
&\textbf{1}& {0} &{1} &{0}&{1}\\
&\textbf{2}& {1}& {1}& {0}&{0}\\
&\textbf{3}&{0}&{0}&{1}&{1}
\end{tabular}
&
\begin{tabular}{@{}cccccc@{}}
&\textbf{B}&&&&\\
\textbf{A}&& \textbf{0} &\textbf{1} &\textbf{2}&\textbf{3}\\
&\textbf{0}& {0}& {1} &{0}&{1}\\
&\textbf{1}& {1} &{0} &{0}&{1}\\
&\textbf{2}& {0}& {1}& {1}&{0}\\
&\textbf{3}&{1}&{0}&{1}&{0}
\end{tabular}&
\begin{tabular}{@{}cccccc@{}}
&\textbf{B}&&&&\\
\textbf{A}&& \textbf{0} &\textbf{1} &\textbf{2}&\textbf{3}\\
&\textbf{0}& {0}& {0} &{1}&{1}\\
&\textbf{1}& {0} &{1} &{1}&{0}\\
&\textbf{2}& {1}& {0}& {0}&{1}\\
&\textbf{3}&{1}&{1}&{0}&{0}
\end{tabular}\\
\end{tabular*}
\end{table}

\section{Discussion}\label{sec6}

We have provided a statistically meaningful interpretation for the building
blocks of GWLP and have generalized resolution by \citet{DenTan99} and \citet{TanDen99} in two meaningful ways for qualitatitve
factors. The complete confounding perspective of $\mathrm{GR}$ of Definition~\ref{de2}
appears to be more sensible than the individual df perspective of
$\mathrm{GR}_{\mathrm{ind}}$ as a primary criterion. However, $\mathrm{GR}_{\mathrm{ind}}$ provides an
interesting new aspect that may provide additional understanding of the
structure of OAs and may help in ranking tied designs. The factor wise
values of Section~\ref{sec5} add useful detail. It will be interesting to pursue
concepts derived from the building blocks of $\mathrm{GR}_{\mathrm{tot}(i)}$ and
$\mathrm{GR}_{\mathrm{ind}(i)}$ for the ranking of mixed level designs. As was
demonstrated in Section~\ref{sec5}, $\mathrm{GR}$ from Definition~\ref{de2} and $\mathrm{GR}_{\mathrm{ind}}$ from
Definition~\ref{de3} are not the only possible generalizations of (\ref{eq4}) for
qualitative factors. The alternative given in equation (\ref{eq6}) appears too
lenient and has therefore not been pursued. The concept of weak strength
deserves further attention: For symmetric designs with weak strength $t$
according to Definition~\ref{de4}, Xu [(\citeyear{Xu03}), Theorem~3] showed that these have
minimum moment aberration (MMA), and consequently GMA (as MMA is equivalent
to GMA for symmetric designs) if they also have maximum $k$-balance for $k
= t+1,\ldots,n$. In particular, this implies that an $\operatorname{OA}(N, s^{n},
t^{-})$ with $N \leq s^{t}$ has GMA, because of Remark~\ref{re4}. Here, we showed
that designs of the highest possible resolution $R$ maximize $\mathrm{GR}$ if they
have weak strength $R$. It is likely that there are further beneficial
consequences from the concept of weak strength.

\begin{appendix}\label{app}
\section*{Appendix: Proof of Theorem~\texorpdfstring{\lowercase{\protect\ref{th1}}}{1}}
Let $\mathbf{M}_{\mathrm{C}} = (\mathbf{1}_{N} ,
\mathbf{M}_{1;\mathrm{C}}, \ldots, \mathbf{M}_{R -1;\mathrm{C}})$, with
$\mathbf{M}_{k;\mathrm{C}}$ the model matrix for all $k$-factor interactions,
$k=1,\ldots,R -1$. The assumption that the resolution of the array is $R$
and the chosen orthogonal contrasts imply $\mathbf{X}_{c}^{\mathrm{T}}
\mathbf{M}_{k;\mathrm{C}} = \mathbf{0}$ for $k < R-1$, with $\mathbf{X}_{c}$
as defined in the theorem. Denoting the $R -1$-factor interaction matrix
$\mathbf{M}_{R -1;\mathrm{C}}$ as $\mathbf{X}_{\mathrm{C}}$, the predictions
for the columns of $\mathbf{X}_{c}$ can be written as
\[
\hat{\mathbf{X}}_{c} = \mathbf{X}_{\mathrm{C}}\bigl (
\mathbf{X}_{\mathrm{C}}^{\mathrm{T}}\mathbf{X}_{\mathrm{C}}  \bigr)^{ -
1}\mathbf{X}_{\mathrm{C}}^{\mathrm{T}}\mathbf{X}_{c} =
\frac{1}{N}\mathbf{X}_{\mathrm{C}}\mathbf{X}_{\mathrm{C}}^{\mathrm{T}}
\mathbf{X}_{c},
\]
since $\mathbf{X}_{\mathrm{C}}^{\mathrm{T}} \mathbf{X}_{\mathrm{C}} = N
\mathbf{I}_{\mathrm{df}(\mathrm{C})}$. As the column averages of $\hat{\mathbf{X}}_{c}$
are 0 because of the coding, the nominators for the $R^{2}$ values are the
diagonal elements of the matrix
\[
\hat{\mathbf{X}}_{c}^{\mathrm{T}}\hat{\mathbf{X}}_{c} =
\frac{1}{N^{2}}\mathbf{X}_{c}^{\mathrm{T}}\mathbf{X}_{\mathrm{C}}
\mathbf{X}_{\mathrm{C}}^{\mathrm{T}}\mathbf{X}_{\mathrm{C}}
\mathbf{X}_{\mathrm{C}}^{\mathrm{T}}\mathbf{X}_{c}\mathop{ =}
_{\mathbf{X}_{\mathrm{C}}^{\mathrm{T}}\mathbf{X}_{\mathrm{C}} =
N\mathbf{I}_{\mathrm{df}(\mathrm{C})}}\frac{1}{N}\mathbf{X}_{c}^{\mathrm{T}}
\mathbf{X}_{\mathrm{C}}\mathbf{X}_{\mathrm{C}}^{\mathrm{T}}
\mathbf{X}_{c}.
\]
Analogously, the corresponding denominators are the diagonal elements of
\[
\mathbf{X}_{c}^{\mathrm{T}}\mathbf{X}_{c} = N
\mathbf{I}_{\mathrm{df}(c)},
\]
which are all identical to $N$. Thus, the sum of the $R^{2}$ values is the
trace of
$\frac{1}{N^{2}}\mathbf{X}_{c}^{\mathrm{T}}\mathbf{X}_{\mathrm{C}}\mathbf{X}_{\mathrm{C}}^{\mathrm{T}}\mathbf{X}_{c}$,
which can be written as
%
\begin{equation}
\operatorname{tr} \biggl( \frac{1}{N^{2}}\mathbf{X}_{c}^{\mathrm{T}}
\mathbf{X}_{\mathrm{C}}\mathbf{X}_{\mathrm{C}}^{\mathrm{T}}
\mathbf{X}_{c} \biggr) = \frac{1}{N^{2}}\operatorname{vec} \bigl(
\mathbf{X}_{\mathrm{C}}^{\mathrm{T}}\mathbf{X}_{c}
\bigr)^{\mathrm{T}}\operatorname{vec} \bigl( \mathbf{X}_{\mathrm{C}}^{\mathrm{T}}
\mathbf{X}_{c} \bigr),\label{eq7}
\end{equation}
where the vec operator stacks the columns of a matrix on top of each other,
that is, generates a column vector from all elements of a matrix [see,
e.g., \citet{Ber09} for the rule connecting trace to vec]. Now, realize
that
\[
\operatorname{vec} \bigl( \mathbf{X}_{\mathrm{C}}^{\mathrm{T}}\mathbf{X}_{c}
\bigr)^{\mathrm{T}} = \operatorname{vec} \Biggl( \Biggl( \sum_{i = 1}^{N}
\mathbf{X}_{\mathrm{C}(i,f)}\mathbf{X}_{c(i,g)} \Biggr)_{(f,g)}
\Biggr)^{\mathrm{T}} = \mathbf{1} {}_{1 \times
N}\mathbf{X}_{u_{1},\ldots,u_{R}},
\]
where an index pair ($i$, $j$) stand for the $i$th row and $j$th column,
respectively, and the columns in\vspace*{1pt} $\mathbf{X}_{u_{1},\ldots,u_{R}}$ are assumed
to appear in the order that corresponds to that in $\operatorname{vec} (
\mathbf{X}_{\mathrm{C}}^{\mathrm{T}}\mathbf{X}_{c}  )^{\mathrm{T}}$
(w.l.o.g.). Then (\ref{eq7}) becomes
\[
\frac{1}{N^{2}}\mathbf{1} {}_{1 \times
N}\mathbf{X}_{u_{1},\ldots,u_{R}}
\mathbf{X}_{u_{1},\ldots,u_{R}}^{\mathrm{T}}\mathbf{1} {}_{1
\times N}^{\mathrm{T}}
= a_{R} ( u_{1},\ldots,u_{R} ),
\]
which proves the assertion.
\end{appendix}



\printaddresses

\end{document}